\documentclass{amsart}

\usepackage{amsmath, amssymb, amsthm}
\usepackage{graphicx}
\usepackage{graphics}
\usepackage{epsfig}
\usepackage[all]{xy}
\usepackage{amsfonts}
\usepackage{hyperref}

\theoremstyle{definition}
\newtheorem{defin}{Definition}
\newtheorem{quest}[defin]{Question}
\newtheorem{prob}[defin]{Problem}
\newtheorem{axi}[defin]{Axiom}
\newtheorem{ex}[defin]{Example}
\newtheorem{rem}[defin]{Remark}
\theoremstyle{plain}
\newtheorem{cor}[defin]{Corollary}
\newtheorem{thm}[defin]{Theorem}
\newtheorem{lem}[defin]{Lemma}
\newtheorem{prop}[defin]{Proposition}

\newcommand{\theorem}[1]{\begin{thm} #1 \end{thm}}

\newcommand{\proposition}[1]{\begin{prop} #1 \end{prop}}

\newcommand{\lemma}[1]{\begin{lem} #1 \end{lem}}
\newcommand{\definition}[1]{\begin{defin} #1 \end{defin}}

\newcommand{\remark}[1]{\begin{rem} #1 \end{rem}}

\newcommand{\corollary}[1]{\begin{cor} #1 \end{cor}}

\newcommand{\lemmaname}[2]{\begin{lem}[#1] #2 \end{lem}}

\newcommand{\ov}[1]{\overline{#1}}

\newcommand{\M}{\mathrm{M}}

\title{Presentations for the Higher Dimensional Thompson Groups $nV$}

\thanks{The authors gratefully acknowledge the partial support by
the NSF grant for Research Experiences for Undergraduates.}

\author{Johanna Hennig}
\address{University of California, San Diego \\ Department of Mathematics, 
 9500 Gilman Drive
La Jolla, CA 92093-0112, USA}
\email{jhennig@math.ucsd.edu}
\author{Francesco Matucci}
\address{University of Virginia, Department of Mathematics,
Charlottesville, VA 22904, USA}
\email{fm6w@virginia.edu}

\begin{document}

\begin{abstract}
In his papers \cite{B1}, \cite{B2} Brin introduced the higher dimensional
Thompson groups $nV$ which are generalizations
to the Thompson group $V$ of self-homeomorphisms of the Cantor set
and found a finite set of generators and relations in the case $n=2$.
We show how to generalize his construction to obtain a finite presentation
for every positive integer $n$. As a corollary, we obtain another proof
that the groups $nV$ are simple (first proved by Brin in \cite{B3}).
\end{abstract}

\maketitle

%\tableofcontents

\section{\label{sec:introduction}Introduction}

The higher dimensional groups $nV$ were introduced by Brin in his papers 
\cite{B1} and \cite{B2} and generalize Thompson's group $V$. We recall that the group $V$ is a group of self-homeomorphisms of the Cantor set $\mathfrak{C}$
that is simple and finitely presented (the standard 
standard introduction to $V$ is the paper by Cannon, Floyd and Parry \cite{CFP}).
The groups $nV$ generalize the group $V$ and act on powers of the Cantor
set $\mathfrak{C}^n$. Brin shows in 
\cite{B1} that the groups $V$ and $2V$ are not isomorphic and
shows in \cite{B2} that the group $2V$ is finitely presented. Bleak and Lanoue \cite{BL}
have recently showed that two groups $mV$ and $nV$ are isomorphic if
and only if $m = n$.

In this paper we give a finite presentation for each of the higher dimensional Thompson groups $nV$.
The argument extends to the ascending union $\omega V$ of the groups $nV$
and returns an infinite presentation of the same flavor. As a corollary, we obtain another proof that the groups
$nV$ and $\omega V$ are simple. Our arguments follow closely and generalize those of Brin in \cite{B1}. \cite{B2} for the group $2V$.

This work arose during a Research Experience for Undergraduates (REU) program at Cornell University. The motivation for the project sprang from a commonly held opinion that the book-keeping required to generalize Brin's presentations to the groups $nV$ 
would be overwhelming.
One would expect from the
similarity of the groups' constructions
 that all arguments for $2V$ would
carry over to $nV$ for all $n$.  Standing in the way of this are the cross
relations.  Thus our paper has two kinds of arguments: those that verify
the parts of \cite{B2} that carry over with no change to $nV$ and those
involving the cross relations that have to be modified to hold in $nV$
(see Lemmas \ref{thm:normalizing-monoid-words} and \ref{thm:brin-lemma-4.22} and Remark \ref{thm:remark-secondary-relation} below).

Following a suggestion of Collin Bleak
the authors have also explored an alternative generating set
(see Section \ref{sec:alternative-generating-set}).
An interesting project would be to find a set of relators for this alternative
generating set in order to use a known procedure which significantly reduces the number of relations,
and which has been successfully 
implemented in a number of papers by Guralnick, Kantor, Kassabov, Lubotzky (see for example \cite{GKKL}).

After a careful reading of Brin's original paper \cite{B2}, it became clear what was needed to generalize his proof, and the current paper borrows heavily from Brin's. Brin was already aware that many of his arguments would probably extend (and he points out in several places in \cite{B1}, \cite{B2}
where it is evident that they do).
We demonstrate how to deal with generators in higher dimensions and what steps are needed to obtain the same type of normalized words which are built for $2V$ in \cite{B2}.

We also mention that Brin asks in \cite{B2} whether or not
the group $2V$ has type $\mathrm{F}_\infty$ (that is, 
having a classifying space that is finite in each dimension).
This has recently been answered by Kochloukova, Martinez-Perez and Nucinkis
\cite{KMN} who have shown that the groups $2V$ and $3V$ have type
$\mathrm{F}_\infty$, therefore obtaining a new proof that these
groups are finitely presented.

\subsection*{Acknowledgments} The authors would like to thank Robert Strichartz  and the National Science Foundation for their support during the REU. The authors would like to thank Collin Bleak and Martin Kassabov for several helpful conversations and
Matt Brin helpful comments and for pointing out that his argument for the simplicity of $2V$ lifts immediately to $nV$ using the presentations
that we find. The authors also would like to thank Matt Brin,
Collin Bleak, Dessislava Kochloukova, Daniel Lanoue, Conchita Martinez-Perez and Brita Nucinkis for kindly referencing the current work while it was still in preparation.
The authors would also like to thank Roman Kogan for advice on how
to create helpful diagrams using Inkscape.

\section{\label{sec:ingredient} The main ingredient and structure of this paper}

Many arguments of Brin generalize word-by-word from $2V$ to $nV$.
For this reason, we advise the reader to have a copy of Brin's papers \cite{B1}, \cite{B2},
as we will adapt some of their results and 
our results will be stated to appear as natural generalizations of those,
including the general argument to show that what we will find is indeed a presentation.

The key observation which allows us to restate many results without proofs (or with little additional effort) is the
following: many statements of Brin do not depend on dimension 2, except
those which need to make use of the ``cross relation'' (relation (18) in Section
\ref{sec:relations} below) to rewrite a
cut in dimension $d$ followed by a cut in dimension $d'$ as one in dimension $d'$ followed by one in dimension $d$.

%We will prove a generalization of this cross-relation which is not an immediate %generalization of Brin's relations for $2V$.
As a result, proofs which need to make use of this new relation require a slight
generalization (for example, the normalization of words in the monoid across
fully divided dimensions) 
while those which do not can be obtained directly using Brin's original
proof. In any case, since statements need to be adapted to our context
we sketch certain proofs to make it clear that they generalize
directly. For example, we will 
show why Brin's proof that $2V$ is simple does not use the new relation (18) 
and therefore it lifts immediately to higher dimensions.

\section{\label{sec:monoid}The monoid $\Pi_n$}

\begin{figure}[0.5\textwidth]
 \begin{center}
  \includegraphics[height=2.3cm]{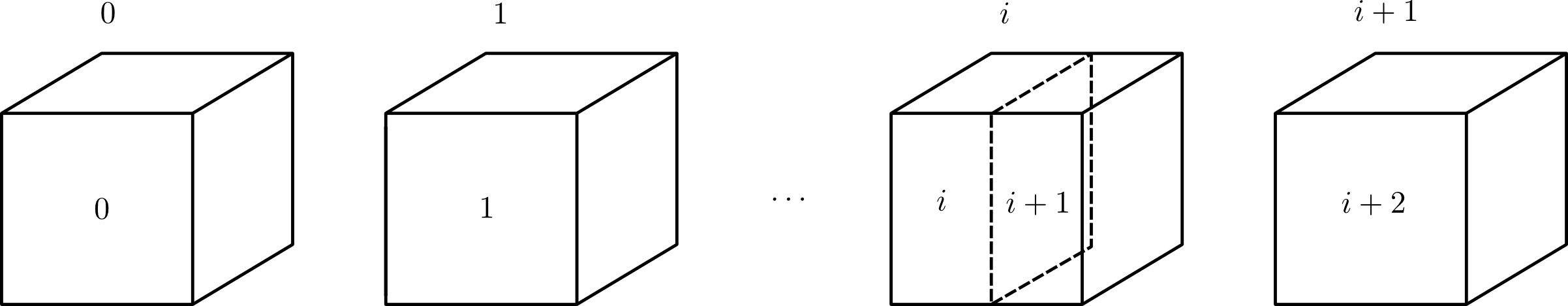}
 \end{center}
\caption{The generator $s_{i,d}$.}
\label{fig:generator-s-i-d}
\end{figure}

\begin{figure}[0.5\textwidth]
 \begin{center}
  \includegraphics[height=2.3cm]{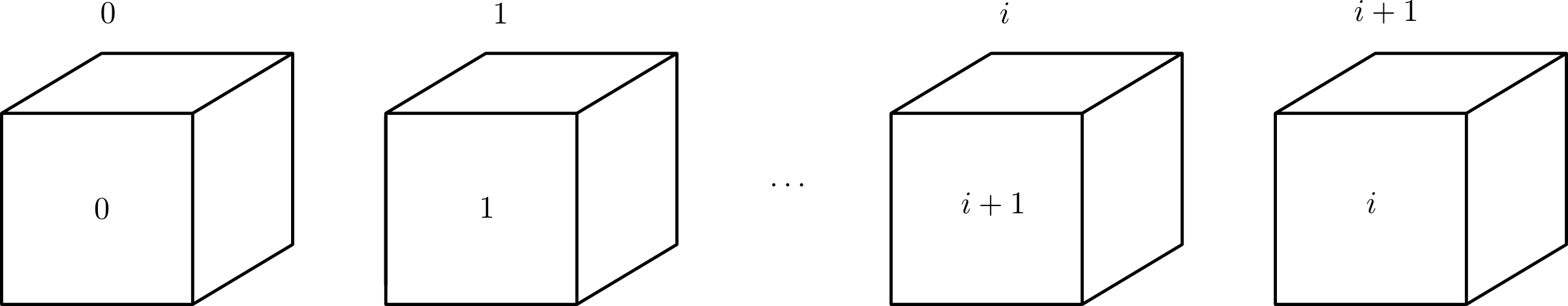}
 \end{center}
\caption{The generator $\sigma_i$.}
\label{fig:generator-sigma-i}
\end{figure}

In \cite{B1} section 4.5, Brin defines the monoid $\Pi$ 
and $\widehat{nV}$ and observes 
that one can extend the definition for all $n$. Elements of $\Pi_n$ are given by numbered patterns in $X$, where $X$ is the union of the set $\{S_0, S_1, . . . \}$ of unit $n$-cubes. Fix $n \in \mathbb{N}$ and fix an ordering on the dimensions $d$, $1 \leq d \leq n$. 
The monoid $\Pi_n$ is generated by the elements $s_{i,d}$ and $\sigma_i$, 
and $s_{i,d}$ denotes the element which cuts the rectangle $S_i$ in half across the $d$-th dimension (see figure \ref{fig:generator-s-i-d})
and $\sigma_i$ is the transposition which switches the rectangle labelled $i$ with that labelled $i+1$, as defined for $2V$ (see figure \ref{fig:generator-sigma-i}).

After each cut, the numbering shifts as before. The following relations hold in $\Pi_n$.
\[
\begin{array}{cccc}
(\M1) & & s_{j,d'}s_{i,d}=s_{i,d}s_{j+1,d'} & i<j, 1 \leq d, d' \leq n \medskip \\

(\M2) & & {\sigma_i}^2 = 1 & i\ge0 \medskip \\

(\M3) & & \sigma_i\sigma_j=\sigma_j\sigma_i & |i-j|\ge2 \medskip \\

(\M4) & & \sigma_i\sigma_{i+1}\sigma_i=\sigma_{i+1}\sigma_i\sigma_{i+1}& i\ge0 \medskip \\

(\M5a) & & \sigma_{j}s_{i,d}=s_{i,d}\sigma_{j+1} & i<j \medskip \\

(\M5b) & & \sigma_{j} s_{i,d}=s_{j+1,d}\sigma_j\sigma_{j+1} & i=j \medskip \\

(\M5c) & & \sigma_{j}s_{i,d}=s_{j,d}\sigma_{j+1}\sigma_j & i=j+1 \medskip \\

(\M5d) & & \sigma_{j}s_{i,d}=s_{i,d}\sigma_j & i>j+1 \medskip \\

(\M6) & & s_{i,d}s_{i+1,d'}s_{i,d'}=s_{i,d'}s_{i+1,d}s_{i,d}\sigma_{i+1} & i\ge0, d \neq d' \\

\end{array}
\]
Note: Relations (M5b) and (M5c) are actually equivalent, using the fact that $\sigma_i$ is its own inverse.

\remark{We observe that the proofs of results of Section 2 in \cite{B2} 
which use relations $(\M1)$ -- $(\M5)$ 
do not depend on the fact that we are in dimension 2,
except for the way they are formulated. For this reason, they 
generalize immediately to the case of the monoid $\Pi_n$
and we do not reprove them. This includes every result up to and including
Lemma 2.9 in \cite{B2}. 

%\marginpar{I'm not sure what the sentence "it requires us to make a choice on how we 
%write elements to the cuts of the ``crossed dimensions''." is trying to say?}
On the other hand, Proposition 2.11 in \cite{B2}
uses the cross relation $(\M6)$ and it requires us to make a choice on how we 
write elements to obtain some underlying pattern. Brin achieves this
type of normalization by writing elements so that vertical cuts appear first, whenever possible. We generalize his argument by describing how to order nodes in forests
(which represent cuts in some dimension).}

The following definition is given inductively on the subtrees.

\definition{Given a forest $F$ we say that a subtree $T$ of some tree of $F$ is \emph{fully divided} across
some dimension $d$ if the root of $T$ is labelled $d$ or if both her left and right subtrees are fully divided across
dimension $d$.}

Given a word $w$ be the word in the generators $\{s_{i,d},\sigma_i\}$,
we define the \emph{length} $\ell(w)$ of $w$ to be the number of appearances in $w$ of elements of
$\{s_{i,d}\}$. It can easily be seen that the length of a word is preserved by relations 
$(\M1)$ -- $(\M6)$.

We restate without proofs Lemmas 2.7, 2.8 and 2.9 from Brin \cite{B2} adapted to our case.

\lemmaname{Brin, \cite{B2}}{\label{thm:brin-lemma-2.7}
If the numbered, labeled forest $F$ comes from a word in $\{s_{i,d} \mid d,i \in
\mathbb{N} \}$, then the leaves of $F$ are numbered so 
that the leaves in $F_i$ have numbers lower
than those in $F_j$ whenever $i < j$ and the leaves in each tree of $F$ are numbered in
increasing order under the natural left right ordering of the leaves.}

\lemmaname{Brin, \cite{B2}}{\label{thm:brin-lemma-2.8}
If two words in the generators $\{s_{i,d}, \sigma_i \mid i \in \mathbb{N},
1 \le d \le n\}$ 
lead to the same numbered, labeled forest, then the words are related by $(\M1)$--$(\M5)$.
}

\lemmaname{Brin, \cite{B2}}{\label{thm:brin-lemma-2.9}
If $F$ is a numbered, labeled forest with the numbering as in Lemma
\ref{thm:brin-lemma-2.7}, and if a linear order is given on 
the interior vertices (and thus of the carets)
of $F$ that respects the ancestor relation, then there is a unique word $w$ in 
$\{s_{i,d} \mid d,i \in \mathbb{N} \}$ leading to $F$ so that the 
order on the interior vertices of $F$ derived from
the order on the entries in $w$ is identical to the given linear order on the interior
vertices.}

The next lemma and corollary are used to prove results analogous to Lemma 2.10 and Proposition 2.11 from \cite{B2}.

\lemma{\label{thm:normalizing-monoid-words}
Let $w$ be a word in the set $\{s_{i,d},\sigma_i\}$ and suppose that the underlying pattern $P$
has a fully divided hypercube $S_i$ across dimension $d$. Then $w \sim w'= s_{i,d}a$ for
some word $a \in \langle s_{i,d},\sigma_i \rangle$.\label{thm:pull-up-divided-dim}}

\begin{proof} 
We use induction on $g:=\ell(w)$. By using relations $(\M5a)$--$(\M5d)$
as in Lemma 2.3 of \cite{B2}
we can assume that $w=pq$ where $p \in \langle s_{i,d} \rangle$ and $q \in \langle \sigma_i\rangle$.
This does not alter the length of $w$. If $g=3$, then $p=p_1p_2p_3$. If $p_1=s_{i,d}$ we are done, otherwise
we have two cases: either $p_2=s_{i+1,d}$ and $p_3=s_{i,d}$ 
or $p_2=s_{i,d}$ and $p_3=s_{i+2,d}$.
Up to using relation $(\M1)$, we can assume that $p_2=s_{i+1,d}$ and $p_3=s_{i,d}$
which is what to want to apply relation $(\M6)$ to $p$ to get $w \sim w'=s_{i,d}s_{i+1,k}s_{i,k}q$.

Now assume the thesis true for all words of length less than $g$. We consider the word $p$ and look at the labelled
unnumbered tree $F_i$ corresponding to $S_i$ with root vertex $u$ and children $u_0$ and $u_1$. Let $T_r$ be the subtree of $F_i$
with root vertex $u_r$, for $r=0,1$. We choose an ordering of the vertices of $F_i$ which respects the ancestor relation and such
that $u$ corresponds to $1$, $u_0$ corresponds to $2$, the other interior nodes of $T_0$ correspond to the numbers
from $3$ to $j=\#($interior nodes of $T_0)$ and $u_2$ corresponds to $j+1$.

By Lemma \ref{thm:brin-lemma-2.9}, the word $p$ is equivalent to
\[
p \sim s_{i,k}(s_{i,m}p_0)(s_{f,l}p_1)
\]
where $s_{i,m}p_0$ is the subword corresponding to the subtree $T_0$ and $s_{f,l}p_1$
is the subword corresponding to the subtree $T_1$ and with $p_0,p_1 \in \langle s_{i,d} \rangle$.
We observe that $\ell(s_{i,m}p_0)<\ell(p)=g$ and $\ell(s_{f,l}\, p_1)<\ell(p)=g$ 
and that
the underlying squares $S_i$ for $s_{i,m}p_0$ and $S_{i+1}$ for $s_{f,l}\,p_1$ are fully divided
across dimension $d$. We can thus apply the induction hypothesis and rewrite
\[
s_{i,m}p_0 \sim s_{i,d}\widetilde{p}_0\widetilde{q}_0 \; \; \; \mbox{and} \; \; \; s_{f,l}\, p_2 \sim s_{f,d}\widetilde{p}_1\widetilde{q}_1.
\]
We restrict our attention to the subword $s_{i,d}\widetilde{p}_0\widetilde{q}_0s_{f,d}$. Using the relations $(\M5a)$--$(\M5d)$
we can move $\widetilde{q}_0$ to the right of $s_{f,d}$ and obtain
\[
s_{i,d}\widetilde{p}_0\widetilde{q}_0s_{f,d} \sim s_{i,d}\widetilde{p}_0s_{g,d}\widetilde{q}
\]
for some permutation word $\widetilde{q}$. 
Since the word $\widetilde{p}_0$ acts on the rectangle $S_i$
and $s_{g,d}$ acts on the rectangle $S_{i+1}$ we can apply
Lemma \ref{thm:brin-lemma-2.8} and \ref{thm:brin-lemma-2.9}
and put a new order on the nodes 
so that the node corresponding to $s_{i,d}$ is $1$
and $s_{g,d}$ is $2$. Thus we have that
\[
s_{i,d}\widetilde{p}_0s_{g,d}\widetilde{q} \sim s_{i,d}s_{i+2,d}\widetilde{p}\, \widetilde{q}
\]
for some $\widetilde{p}$ word in the set $\{s_{i,d}\}$. Thus we have
$w \sim w''=s_{i,k}s_{i,d}s_{i+2,d}\widetilde{p}\, \widetilde{q}$ and so, by applying the cross relation $(\M6)$
to the first three letters of $w''$ we get
\[
w \sim w'' \sim w'=s_{i,d}s_{i,k}s_{i+2,k}\widetilde{p}\, \widetilde{q}=s_{i,d}a \qquad 
\]
\end{proof}

We have now proved Lemma 2.10 from \cite{B2}, since in order for a tree in a forest to be non-normalized, one of the rectangles in the pattern corresponding to that tree must be fully divided across two different dimensions.

\lemmaname{Brin, \cite{B2}}{\label{thm:brin-lemma-2.10}If 
two different forests correspond to the same pattern in $X$, then
at least one of the two forests is not normalized.}

\remark{Lemma \ref{thm:normalizing-monoid-words} is used in our extension of Brin 2 Proposition 2.11 so that we can push dimension $d$ under the root.
This is explained better in the following Corollary.}

\corollary{Let $w$ be a word in the generators $\{s_{i,d},\sigma_i\}$ such that its underlying square $S_i$
is fully divided across dimensions $d$ and $\ell$. Then
\[
w \sim w'=s_{i,d}s_{i,\ell}s_{i+2,\ell}a \sim w'' = s_{i,\ell}s_{i,d}s_{i+2,d}b
\]
for some suitable words $a,b$ in the generators $\{s_{i,d},\sigma_i\}$.\label{thm:pull-up-fully-divided}}

\begin{proof} This is achieved by a repeated application of the previous Lemma \ref{thm:normalizing-monoid-words}. 
We apply Lemma \ref{thm:normalizing-monoid-words}
to $w$ and obtain $w\sim s_{i,d}a_1$. By construction, we notice that 
the underlying squares $S_i$ and $S_{i+1}$ of $a_1$
are fully divided across dimension $\ell$, so we can apply the previous Lemma to $a_1$ to get $a_1 \sim s_{i,\ell}a_2$
and finally we apply it again to $a_2 \sim s_{i+2,\ell}a$. Hence $w \sim w'=s_{i,\ell}s_{i+2,\ell}a$.
To get $w''$ we apply the cross relation $(\M6)$ 
to the subword $s_{i,\ell}s_{i,d}s_{i+2,d}$. 
\end{proof}

\proposition{\label{thm:analogue-brin-prop-2.11}A word w is related by 
$(\M1)$ through $(\M6)$ to a word corresponding to a normalized, labelled forest.}

\begin{proof} We proceed by induction on the length of $w$. Let $g$ be the length of $w$
and assume the result holds for all words of length less than $g$. 
As before, write $w=pq$, 
where $p=s_{i_0}s_{i_1}\dots s_{i_{n-1}}$ (here, the $i_j$ refers to the cube which is being cut; 
we omit the second index indicating dimension as it is unimportant for now). Write $w=s_{i_0}w'$; 
since the order of the interior vertices of the forest for $p$ given by the order of the letters in $p$ 
must respect the ancestor relation, we know that the interior vertex corresponding to $s_{i_0}$ must be a 
root of some tree, $T$. As $w'$ is a word of length less than $g$, we may apply our inductive hypothesis 
and assume that $w'$ can be rewritten via relations $(\M1)$ through 
$(\M6)$ to obtain a 
corresponding normalized forest. The pattern $P$
for \(w\) is obtained from the pattern $P'$ for $w'$ by
applying the pattern of $P'$ in unit square $S_i$ to the
rectangle numbered $i$ in the pattern for $s_{i_0}$.  The forest $F$
for \(w\) is obtained from the forest $F'$ for $w'$ by
attaching the $i$-th tree of $F'$ to the $i$-th leaf of the
forest for $s_{i_0}$.  Since $F'$ is normalized, it is seen that $F$
has all interior vertices normalized except possibly for the root
vertex of one tree, $T$.

Let $u$ be the root vertex of $T$ with label $k$ and with 
children $u_1$ and $u_2$. Let $T_1$ and $T_2$ be the subtrees of $T$ whose roots 
are $u_1$ and $u_2$, respectively. By hypothesis, $T_1$ and $T_2$ are already normalized. 
If $T$ is not normalized already, then $T$ must be fully divided across the dimension that $u$ is 
labeled with, $k$, and some other dimension less than $k$. Let $d$ be the minimal dimension across 
which $T$ is fully divided. Since $T_1$ and $T_2$ are also fully divided across $d$, by Lemma \ref{thm:pull-up-divided-dim}, 
we may apply relations $(\M1)$ through $(\M6)$ 
to the subwords of $w$ corresponding to $T_1$, $T_2$ until $u_1$ 
and $u_2$ are each labelled $d$. Now by lemma 2.9, we may assume $w= s_{i_0, k}s_{i_0, d}s_{i_0+2, d}w''$ 
where $w''$ is the remainder of $w$. We apply relation $(\M6)$ to obtain $w= s_{i_0, d}s_{i_0, k}s_{i_0+2, k}\sigma_{i_0}w''$. 
Now, we have normalized the vertex $u$, and we may now use the inductive hypothesis to renormalize the trees $T_1$ 
and $T_2$. The result is a normalized forest.
\end{proof}

The proof of the following result follows the same argument of Theorem 1
in \cite{B2}, using Lemma 2.10 in \cite{B2} and Proposition \ref{thm:analogue-brin-prop-2.11} (to extend Proposition 2.11 in \cite{B2}).

\theorem{The monoid $\Pi_n$ is presented by using the generators 
$\{s_{i,d}, \sigma_i\}$ and relations $(\M1)$--$(\M6)$.}

\section{\label{sec:relations}Relations in $nV$}

\subsection{Generators for $nV$}

The following generators are defined as in \cite{B1} and analogous arguments show why they are a generating set for $nV$.

\[
\begin{array}{ccc}
X_{i,d} = (s_{0,1}^{i+1} s_{1,d}, s_{0,1}^{i+2}) & i \ge 0, 1 \le d \le n & {}\bigskip \medskip \\

C_{i,d} = (s_{0,1}^i s_{0,d}, s_{0,1}^{i+1}) & i \ge 0, 2 \le d \le n & (\mbox{baker's maps}) \bigskip \medskip \\ 

\pi_i = (s_{0,1}^{i+2} \sigma_1, s_{0,1}^{i+2}) & i \ge 0 \; \; \; & (\mbox{$\sigma_i$ defined as above}) \bigskip \medskip \\

\overline{\pi}_i = (s_{0,1}^{i+1} \sigma_0, s_{0,1}^{i+1}) & i \ge 0 & {}
\end{array}
\]

\subsection{Relations involving cuts and permutations}
In all the following relations (1) -- (7) the reader can assume that 
$1 \le d, d' \le n$, unless otherwise stated.
\[
\begin{array}{cccc}
(1) & & X_{q,d} X_{m,d'} = X_{m,d'} X_{q+1,d} & m < q, \medskip \\

(2) & & \pi_q X_{m,d} = X_{m,d} \pi_{q+1} & m < q \medskip \\

(3) & & \pi_q X_{q,d} = X_{q+1,d} \pi_q \pi_{q+1} & q \ge 0 \medskip \\

(4) & & \pi_q X_{m,d} = X_{m,d} \pi_q  & m > q+1 \medskip \\

(5) & & \overline{\pi}_q X_{m,d} = X_{m.d} \overline{\pi}_{q+1} & m < q \medskip \\

(6) & &  \overline{\pi}_m X_{m,1} = \pi_m \overline{\pi}_{m+1} & m \ge 0 \medskip \\

(7) & & X_{m,d} X_{m+1,d'} X_{m,d'} = X_{m,d'} X_{m+1,d} X_{m,d} \pi_{m+1} & m \ge 0, d \ne d'  \medskip \\ 

\end{array}
\]

\subsection{Relations involving permutations only}

\[
\begin{array}{cccc}
(8) & & \pi_q \pi_m = \pi_m \pi_q & |m-q| > 2 \medskip \\

(9) & & \pi_m \pi_{m+1} \pi_m = \pi_{m+1} \pi_m \pi_{m+1} & m \ge 0 \medskip \\

(10) & & \overline{\pi}_q \pi_m = \pi_m \overline{\pi}_q & q \ge m+2 \medskip \\

(11) & & \pi_m \overline{\pi}_{m+1} \pi_m = \overline{\pi}_{m+1} \pi_m \overline{\pi}_{m+1}  & m \ge 0 \medskip \\

(12) & & \pi_m^2 =1 & m \ge 0 \medskip \\

(13) & & \overline{\pi}_m^2 =1 & m \ge 0 \medskip \\

\end{array}
\]

\subsection{Relations involving baker's maps}
In all the following relations (14) -- (18) the reader can assume that 
$2 \le d \le n$ and $1 \le  d' \le n$, unless otherwise stated.
\[
\begin{array}{cccc}

(14) & &  \overline{\pi}_m X_{m,d} = C_{m+1,d} \pi_m \overline{\pi}_{m+1}  & m \ge 0,  \medskip \\

(15) & & C_{q,d} X_{m,d'} = X_{m,d'} C_{q+1,d} & m < q, \medskip \\

(16) & & C_{m,d} X_{m,1} = X_{m,d} C_{m+2,d} \pi_{m+1} & m \ge 0, \medskip \\

(17) & & \pi_q C_{m,d} = C_{m,d} \pi_q & m > q+1 \medskip \\

(18) & & C_{m,d} X_{m,d'} C_{m+2,d'}= C_{m,d'} X_{m,d} C_{m+2,d} \pi_{m+1} & m \ge 0, 1 < d' < d \le n \medskip \\
\end{array}
\]

Relations (1) through (17) are generalizations of those given in \cite{B1} and their proofs are completely analogous. The only new family of relations is (18) which we prove using relation $(\M6)$ from the monoid:

\begin{proof}

\[
\begin{split}
C_{m,d} X_{m,d'} C_{m+2,d'} &= (s_{0,1}^ms_{0,d},\,
s_{0,1}^{m+1}) (s_{0,1}^{m+1}s_{1,d'},\, s_{0,1}^{m+2}) (s_{0,1}^{m+2}s_{0,d'},\, s_{0,1}^{m+3}) \\ &= (s_{0,1}^ms_{0,d}s_{1,d'}s_{0,d'},\,
s_{0,1}^{m+3}) \\ &= (s_{01}^ms_{0,d'}s_{1,d}s_{0,d}\sigma_1,\, s_{0,1}^{m+3}) \\ &=
 (s_{0,1}^ms_{0,d'},\,
s_{0,1}^{m+1}) (s_{0,1}^{m+1}s_{1,d},\, s_{0,1}^{m+2}) (s_{0,1}^{m+2}s_{0,d},\, s_{0,1}^{m+3}) (s_{0,1}^{m+3}\sigma_1,\, s_{0,1}^{m+3}) \\ &= C_{m,d'} X_{m,d} C_{m+2,d} \pi_{m+1}. 
\end{split}
\]
\end{proof}

\lemmaname{Subscript Raising Formulas}{\label{thm:subscript-raising}We have that
\[
\begin{array}{c}
C_{r,d} \sim C_{r+1, d} X_{r,d} \pi_{r+1}X_{r, 1}^{-1} \\
{} \\
\overline{\pi}_r \sim \pi_r \overline{\pi}_{r+1}X_{r,1}^{-1} \sim X_{r,1}\overline{\pi}_{r+1}\pi_r 
\end{array}
\]}

We observe that the first formula of Lemma \ref{thm:subscript-raising} follows from relations (15) and (16), while the second is a generalization of the one found in \cite{B2}.

\subsection{Secondary Relations for $nV$}

\[
\begin{split}
%%%%%%%%%%%%%%%%%%%%
%%%%%%%%%%%%%%%%%%%%
%%%%%%%%%%%%%%%%%%%%
X_{q,d}^{-1}X_{r,d} & \sim
\begin{cases}
X_d X_d^{-1} & r\ne q \\
1  & r=q
\end{cases}
\qquad (1 \le d \le n) \\
%%%%%%%%%%%%%%%%%%%%
%%%%%%%%%%%%%%%%%%%%
%%%%%%%%%%%%%%%%%%%%
X_{q,d}^{-1}X_{r,d'} & \sim 
\begin{cases}
X_{d'} X_d^{-1} & r \ne q %\qquad(\mbox{uses (1)})
\\
w(X_{d'})\pi w(X_d^{-1})  & r = q %\qquad (\mbox{uses (7)})
\end{cases} 
\qquad (1\le d, d' \le n, d \ne d') \\
%%%%%%%%%%%%%%%%%%%%
%%%%%%%%%%%%%%%%%%%%
%%%%%%%%%%%%%%%%%%%%
%C_{q,d}^{-1}X_{r,1} & \sim 
%\begin{cases}
%X_1 C_d^{-1} & r < q \\%\qquad(\mbox{uses (15)})\\
%w(X_1,\pi,X_d^{-1})C_d^{-1}  & r \ge q %\qquad (\mbox{subscript raising formula, then (15)})
%\end{cases} \\
%%%%%%%%%%%%%%%%%%%%
%%%%%%%%%%%%%%%%%%%%
%%%%%%%%%%%%%%%%%%%%
C_{q,d}^{-1}X_{r,d'} & \sim 
\begin{cases}
X_{d'} C_d^{-1} & r < q \\
w(X_1, \pi,X_d^{-1}) X_{d'} C_d^{-1}  & r \ge q %\qquad (\mbox{check this})
\end{cases}
\qquad (2 \le d \le n, 1 \le d' \le n) \\
%%%%%%%%%%%%%%%%%%%%
%%%%%%%%%%%%%%%%%%%%
%%%%%%%%%%%%%%%%%%%%
%X_{r,1}^{-1}C_{q,d} & \sim
%\begin{cases}
%C_d X_1^{-1} & r<q \\
%C_d w(X_d,\pi,X_1^{-1}) & r \ge q
%\end{cases} \\
%%%%%%%%%%%%%%%%%%%%
%%%%%%%%%%%%%%%%%%%%
%%%%%%%%%%%%%%%%%%%%
X_{r,d'}^{-1}C_{q,d} & \sim
\begin{cases}
C_d X_{d'}^{-1} & r<q \\
C_d X_{d'}^{-1} w(X_{d},\pi,X_1^{-1}) & r \ge q
\end{cases}
\qquad (2 \le d \le n, 1 \le d' \le n) \\
%%%%%%%%%%%%%%%%%%%%
%%%%%%%%%%%%%%%%%%%%
%%%%%%%%%%%%%%%%%%%%
%C_{q,i}^{-1}X_{r,j} & \sim 
%\begin{cases}
%X_j C_i^{-1} & r < q \\
%w(X_1,\pi,X_i^{-1})w(X_i) C_i^{-1}  & r \ge q
%\end{cases} 
%\qquad (i \ne j \ne 1) \\
%%%%%%%%%%%%%%%%%%%%
%%%%%%%%%%%%%%%%%%%%
%%%%%%%%%%%%%%%%%%%%
\pi_q X_{r,d} & \sim \begin{cases}
X_d w(\pi) & \qquad (1 \le d \le n)
\end{cases} \\
%%%%%%%%%%%%%%%%%%%%
%%%%%%%%%%%%%%%%%%%%
%%%%%%%%%%%%%%%%%%%%
\overline{\pi}_q X_{r,1} & \sim
\begin{cases}
X_1 \overline{\pi} & r < q \\
\pi \overline{\pi} & r=q \\
w(X_1)\overline{\pi}w(\pi) & r>q
\end{cases} \\
%\end{split}
%\]
%
%\[
%\begin{split}
%%%%%%%%%%%%%%%%%%%%
%%%%%%%%%%%%%%%%%%%%
%%%%%%%%%%%%%%%%%%%%
\overline{\pi}_q X_{r,d} & \sim
\begin{cases}
X_d \overline{\pi} & r < q \\
C_d \pi \overline{\pi} & r=q \\
w(X_1)X_d\overline{\pi}w(\pi) & r>q
\end{cases}
\qquad (2 \le d \le n) \\
%%%%%%%%%%%%%%%%%%%%
%%%%%%%%%%%%%%%%%%%%
%%%%%%%%%%%%%%%%%%%%
\pi_q C_{r,d} & \sim 
\begin{cases}
C_d \pi & r>q+1 \\
C_d w(X_1^{-1}, \pi, X_d) & r \le q+1
\end{cases} 
\qquad (2 \le d \le n) \\
%%%%%%%%%%%%%%%%%%%%
%%%%%%%%%%%%%%%%%%%%
%%%%%%%%%%%%%%%%%%%%
\overline{\pi}_q C_{r,d} & \sim
\begin{cases}
X_d \overline{\pi}\pi & r=q+1 \\
w(X_1)X_d \overline{\pi}w(\pi) & r>q+1 \\
w(X_d) C_d \pi \overline{\pi} w(\pi,X_1^{-1}) & r<q+1
\end{cases} 
\qquad (2 \le d \le n) \\
%%%%%%%%%%%%%%%%%%%%
%%%%%%%%%%%%%%%%%%%%
%%%%%%%%%%%%%%%%%%%%
%X_{q,i}^{-1}C_{r,j} & \sim
%\begin{cases}
%C_j X_i^{-1} & q<r \\
%C_j X_i^{-1} w(X_j,\pi,X_1^{-1}) & q>r
%\end{cases} 
%\qquad (1 \ne i \ne j) \\
%%%%%%%%%%%%%%%%%%%%
%%%%%%%%%%%%%%%%%%%%
%%%%%%%%%%%%%%%%%%%%
C_{q,d}^{-1}C_{r,d} & \sim
\begin{cases}
w(X_1^{-1},\pi,X_d) & q<r \\
1 & q=r \\
w(X_1,\pi,X_d^{-1}) & q>r
\end{cases} 
\qquad (2 \le d \le n) \\
%%%%%%%%%%%%%%%%%%%%
%%%%%%%%%%%%%%%%%%%%
%%%%%%%%%%%%%%%%%%%%
C_{q,d}^{-1}C_{r,d'} & \sim
\begin{cases}
X_{d'} C_{d'} \pi C_d^{-1}X_{d}^{-1} w(X_{d'}, \pi, X_{1}^{-1}) & q>r \\
X_{d'}C_{d'}\pi C_{d}^{-1}X_{d}^{-1} & q=r \\
w(X_{1}, \pi, X_{d'}^{-1}) X_{d}C_{d}\pi C_{d'}^{-1}X_{d'}^{-1} & q<r
\end{cases}
\qquad (1 \le d' < d \le n)
\end{split}
\]

%\marginpar{New secondary relation, new proof and new remark 12}
\begin{proof}
We only prove the last set of secondary relations as it is the only one that
does not immediately descend from the computations in Brin \cite{B2}.
If $q>r$ we can apply the subscript raising formulas repeatedly for
$j$ times until $r+j=q$ and rewrite the product as
\[
C_{q,d}^{-1}C_{r,d'} \sim C_{q,d}^{-1} C_{r+1,d'} X_{r,d'} \pi_{r+1} X_{r,1}^{-1}
\sim \ldots \sim C_{q,d'}^{-1} C_{r+j,d'} w(X_{d'}, \pi, X_{1}^{-1}).
\]
We argue similarly if $q<r$. We now have to study the product
$C_{q,d}^{-1}C_{q,d'}$. Without loss of generality we assume $d' < d$
and apply relation (18):
\[
C_{q,d}^{-1}C_{q,d'} = X_{q,d'}C_{q+2,d'}\pi_{q+1}C_{q+2,d}^{-1}X_{q,d}^{-1},
\]
which is what was claimed. Similar relations can be derived if $d' > d$.
\end{proof}

\remark{\label{thm:remark-secondary-relation}When using the the 
last two secondary relations, we alter a word in a way that does not increase
the number of $C$'s. This allows us to generalize the proof of Lemma 4.6 in Brin
\cite{B2} thus rewriting a word
of type $w(X, C, \pi, C^{-1}, X^{-1})$ in $LMR$ form 
so that the number of $C$'s does not increase
(see Lemma \ref{thm:brin-lemma-4.6} below).
This observation lets us generalize Lemma 4.7 in Brin \cite{B2}
(see Lemma \ref{thm:brin-lemma-4.7} below).
%the letters $\ov{\pi}$ appear in the word $w$ and one needs to take $w$
%into $LMR$ form again. 
In fact, all our secondary relations are immediate
generalizations of those in Brin \cite{B2} and the last one
does not introduce appearances of $\ov{\pi}$ and therefore all
the letters in the last secondary relations can be migrated to their needed
position by means of the previous secondary relations, without altering
the original argument of Lemma 4.7 in Brin \cite{B2}. 
Therefore even in the case of $nV$ one is able to do the book-keeping
without risk of creating extra letters which cannot be passed safely without
recreating them, and hence we obtain an argument which terminates.
}

\section{\label{sec:presentation-proof}Presentations for $nV$}

We now show how the relations above enable us to put our group elements into a normal form, starting with words in the generators of $nV$ corresponding to elements from $\widehat{nV}$.

\lemma{\label{thm:brin-lemma-4.1}Let $w$ be a word in $\{X_{i,d}, \pi_i, 
X_{i,d}^{-1} 
\mid 1 \leq d \leq n, i \in \mathbb{N} \}$. Then $w \sim LMR$ where $L$ and $R^{-1}$ are words in $\{X_{i,d}\}$ and $M$ is a word in $\{ \pi_i \}$.}

\begin{proof} There is a homomorphism from $\widehat{nV}$ to $nV$ given by $s_{i,d} \mapsto X_{i,d}$ and $\sigma_i \mapsto \pi_i$. This follows from the correspondence between the relations for $\widehat{nV}$ and $nV$ as given below:
\[
\begin{split} (\M1)&\rightarrow (1), \\
(\M2) &\rightarrow (12), \\ (\M3)
&\rightarrow (8), \\ (\M4) &\rightarrow (9),
\end{split}\qquad\qquad\qquad
\begin{split} (\M5a) &\rightarrow (2),
\\ (\M5b), (\M5c) &\rightarrow (3), \\
(\M5d) &\rightarrow (4), \\
(\M6) &\rightarrow (7).  
\end{split}
\]
Hence, any word $w$ as given above is the image under this homomorphism of a word $w'$ in $\widehat{nV}$. Since $\widehat{nV}$ is the group of right fractions of the monoid $\Pi_n$, we can represent $w'$ as $pq^{-1}$ where $p$, $q$ are words in $\{s_{i,d}, \sigma_i \mid 1 \leq d \leq n, i \in \mathbb{N} \}$. Now, as noted before in the proof of Lemma \ref{thm:normalizing-monoid-words}, we can assume $p$ and $q$ are in the form $ab$ where $a \in \langle s_{i,d} \rangle$ and $b \in \langle \sigma_i \rangle$. Hence, we have written $w'$ as $lmr$ for $l, r^{-1}$ $\in \langle s_{i,d} \rangle$ and $m \in \langle \sigma_i \rangle$ since elements of $\langle \sigma_i \rangle$ are their own inverse. Applying the homomorphism to $w'$ puts $w$ in the desired form. 
\end{proof}

The following two results follow the original proofs of Lemma 4.6 and 4.7
in Brin \cite{B2} via Remark \ref{thm:remark-secondary-relation}.

\lemma{\label{thm:brin-lemma-4.6}Let 
$w$ be of the form $w(X, C, \pi, X^{-1}, C^{-1})$. 
Then $w \sim LMR$ where $L$ and $R^{-1}$ 
are words of the form $w(X, C)$ and $M$ is of the form $w(\pi)$.
Further the number of appearances of $C$ in $L$ will be no larger than the number of
appearances of $C$ in $w$ and the number of appearances of $C^{-1}$ in $R$ 
will be no larger than the number of appearances of $C^{-1}$ in $w$.
}

\lemma{\label{thm:brin-lemma-4.7}Let 
$w$ be a word in the generating set
$\{X_{i,d}, C_{i,d'}, \pi_i, \ov{\pi}_i,, X_{i,d}^{-1}, C_{i,d'}^{-1} 
\mid 1 \leq d \leq n, 2\le d' \le n, 
i \in \mathbb{N} \}$. Then $w \sim LMR$ where $L$ and
$R^{-1}$ are words of the form $w(X, C)$ and $M$ is of the form 
$w(\pi, \ov{\pi})$.}

\lemma{\label{thm:brin-lemma-4.8}Let 
$w$ be a word in the generating set 
\[
\{X_{i,d}, C_{i,d'}, \pi_i, \ov{\pi}_i,, X_{i,d}^{-1}, C_{i,d'}^{-1} 
\mid 1 \leq d \leq n, 2\le d' \le n, 
i \in \mathbb{N} \}.
\] 
Then $w \sim LMR$ where  
\begin{itemize}
\item $L=C_{i_0, d_0}C_{i_1, d_1}\dots C_{i_g, d_g}q$ with $i_0<i_1<\cdots<i_g$ for $g \geq -1$ and $q$ is a word in the set $\{X_{i,d} \mid 1 \leq d \leq n, i \in \mathbb{N}\}$
\medskip
\item $R^{-1}=C_{j_0, d'_0}C_{j_1, d'_1}\dots C_{j_m, d'_m}q'$ with $j_0<j_1<\cdots<j_m$ for $m \geq -1$ and $q'$ is a word in the set $\{X_{i,d} \mid 1 \leq d \leq n, i \in \mathbb{N}\}$
\medskip
\item M is a word in the set $\{ \pi_i, \overline{\pi}_i \mid i \in \mathbb{N} \}$
\end{itemize}
}

\begin{proof}
By using the secondary relations, we can assume that $w \sim LMR$ where $L$ and $R^{-1}$ are words in $\{X_{i,d}, C_{i,d}\}$ and $M$ is a word in $\{ \pi_i, \overline{\pi}_i\}$ by analogous arguments used in 
Lemmas 4.6 and 4.7 of \cite{B2}. We then improve $L$ using the subscript raising formula for the $C_{i,d}$ and relation (15) as in the proof of Lemma 4.8 of \cite{B2}.
We notice that to adapt the quoted lemmas from \cite{B2} we need to make
use of Remark \ref{thm:remark-secondary-relation} to make sure that
the appearances of $C$'s and $\overline{\pi}$'s do not increase.
\end{proof}

We define the notions of \emph{primary} and \emph{secondary tree} 
and of \emph{trunk} exactly
the same way that Brin does in \cite{B2}.
The primary tree is the tree corresponding to the word $t$ in Lemma 18
and any extension to the left is a secondary tree for $L$.
The following extends Lemma 4.15 \cite{B2} adapted to our case. The proof is completely analogous.

\lemma{\label{thm:brin-lemma-4.15}Let $L=C_{i_0, d_0}C_{i_1, d_1}\cdots
C_{i_g, d_g} X_{i_{n+1}, d_{n+1}} \cdots X_{i_{l-1}, d_{l-1}}$ 
where \(i_0<i_1<\cdots < i_g\), where $2 \leq d_k \leq n$
for $k \in \{0,\dots, g\}$ and $1 \leq d_k \leq n$ for $k \in \{g+1,\dots, l-1\}$.  Let
\(m\) equal the maximum of 
\[
\{i_j+g+2-j\mid g+1\le j\le l-1\}\cup\{i_g+1\}.
\] 
Then \(L\) can be represented as \(L=(t,
s_{0,1}^k)\) where \(t\) is a word in $\{s_{i,d} \}$ and \(k\) is the length
of \(t\), so that \(k=m+l-g\), and so that the tree \(T\) for \(t\)
is the primary tree for \(L\) and is described as follows.  The tree
\(T\) consists of a trunk \(\Lambda \) with a finite forest \(F\)
attached.  The trunk \(\Lambda \) has \(m\) carets and \(m+1\)
leaves numbered \(0\) through \(m\) in the right-left order.  If the
carets in \(\Lambda \) are numbered from 0 starting at the top, then
the label of the \(i\)-th caret is \(d_k\) if \(i = i_k\) for k in \(\{0,1,
\ldots\, g\}\) and \(1\) otherwise.} 

The following two lemmas are used in proving 
Proposition \ref{thm:remark-secondary-relation} which allows us to assume the trees corresponding to our group elements are in normal form.  

\lemma{\label{thm:brin-lemma-4.21}Let \(L=C_{i_0, d_0}C_{i_1, d_1}\cdots
C_{i_g, d_g}u\) and \(L'=C_{k_0, d'_0}C_{k_1, d'_1}\cdots C_{k_g, d'_g}u'\) where
\(i_0<i_1<\cdots<i_g\), where \(k_0<k_1<\cdots <k_g\), where \(u\)
is a word in the set $\{X_{i,d} \mid 1 \leq d \leq n, i \in \mathbb{N}\}$, and where \(u'\) is a word in the set
$\{X_{i,d}, \pi_i \mid 1 \leq d \leq n, i \in \mathbb{N}\}$.  Assume that \(L\) is expressible as \((t,s_{0,1}^p)\)
as an element of \(\widehat{nV}\) with \(t\) a word in
$\{s_{i,d}\}$ and \(p\) is the length of \(t\).  Let \(m\) be the
number of carets of the trunk of the tree \(T\) corresponding to $t$ and assume that \(m\ge
k_g+1\).

If \(L\sim L'\), then there is a word \(u''\) in $\{X_{i,d}\}$,
and there is a word \(z\) in \(\{\pi_i\mid i\le p-2\}\) so that
setting \(L_1=C_{k_0, d'_0}C_{k_1, d'_1}\cdots C_{k_g, d'_g}u''\) and \(L_2=L_1z\)
gives that \(L\sim L_2\) and \(L_1\) is expressible as \((t',
s_{0,1}^p)\) with \(t'\) a word in $\{s_{i,d}\}$ of length \(p\) so
that the tree \(T'\) for \(t'\) is normalized except possibly at
interior vertices in the trunk of the tree, and so that the trunk of
\(T'\) has \(m\) carets. }

\begin{proof} The homomorphism $\widehat{nV} \rightarrow nV$ given by $s_{i,d} \mapsto X_{i,d}$ and $\sigma_i \mapsto \pi_i$ allows us to write $u' \sim u''z'$ with $u''$ a word in $\{X_{i,d}\}$ and $z'$ is a word in \(\{\pi_i\mid i \in \mathbb{N}\}\) such that the forest $F$ for $u''$ is normalized. The rest of the proof goes through as before, but we describe the slight modifications needed for our case. We write $L = (t s_{0,1}^k,s_{0,1}^{p+k} ) = (\widehat{t}s_{1,0}^r x, s_{1,0}^{q+r}) = L_2$ as elements in $\widehat{nV}$ where $x$ is a word in $\{ \sigma_i \}$ and $p +k = q+r$. As before, we can conclude that the unnumbered patterns for $t s_{0,1}^k$ and $\widehat{t}s_{1,0}^r$ are identical. 

In the tree for $t s_{0,1}^k$, let the left edge vertices be $a_0, a_1, \dots, a_b$ reading from the top, so that $a_0$ is the root of the tree. Since we assume the trunk of the tree has $m$ carets, we know $b=m+k$  and for $m \leq i < b$, the label for $a_i$ is $1$. Similarly, in the tree for $\widehat{t}s_{1,0}^r$, let the left  edge vertices be $a_0', a_1', \dots, a_b'$ reading from the top. Note that remark (*) in the proof
of Theorem 4.21 in Brin \cite{B2} (which we are about to restate)
remains true in our general case, by giving a new definition: for each left edge vertex, $a_i$, define the $n$-tuple $(x_1^i, \dots, x_n^i)$ where $x_k^i$ equals the number of left edge vertices above $a_i$ with label $k$. (Note we are using $i$ to denote an index, not an exponent). It follows that $x_1^i + x_2^i + \dots + x_n^i$ is the total number of left edge vertices above $a_i$. Then we have:

\medskip
(*) The rectangle corresponding to a left edge vertex $a_i$ depends only on the $n$-tuple $(x_1^i, \dots, x_n^i)$
\medskip

In other words, for the rectangle labeled $''0''$ in any pattern, the order of the different cuts does not matter. This is because the rectangle labeled $''0''$ must contain the origin and its size in each dimension $k$ will be $2^{-x_k^i}$. Hence, the analogous statement for our case follows, and we conclude that the $n$-rectangle $R$ corresponding to $a_m$ is identical to the $n$-rectangle $R'$ corresponding to $a_m'$ Since $R$ is divided $k$ times across dimension $1$, so is $R'$, and hence the tree below $a_m'$ must consist of an extension to the left by $k$ carets all labeled $1$, and we can conclude that $r \geq k$. The rest of the proof follows exactly as before. 
\end{proof}

Here, we define a notion of \emph{complexity} 
to measure progress in the following lemma and proposition towards normalizing trees. If $T$ is a labeled tree, let $a_0, a_1, \dots, a_m$ be the interior, left edge vertices of $T$ reading from top to bottom so that $a_0$ is the root. Let $b_0b_1 \dots b_m$ be a word in $\{ 1, 2, \dots, n \}$ where $b_i = k$ if $a_i$ is labelled $k$ for $0 \leq i \leq m$. We say $b_0b_1 \dots b_m$ is the complexity of $T$. We impose the length-lex ordering
on such words, that is if $w_1$ and $w_2$ are two such words, then we
say $w_1 < w_2$ if $w_1$ is shorter than $w_2$ or
if $w_1 = b_0^1 \dots b_m^1$ and $w_2 = b_0^2 \dots b_m^2$ are two such words of the same length, then $w_1 < w_2$ if when we take $j \in \{0, \dots, m \}$ minimal where $b_j^1 \neq b_j^2$, we have $b_j^1 < b_j^2$. We will refer to this notion in the following lemma.

\lemma{\label{thm:brin-lemma-4.22}Let \(L=C_{i_0,d_0}C_{i_1,d_1}\cdots
C_{i_g,d_g}u\) where \(i_0<i_1<\cdots <i_g\) and $u$ is a word in the set 
$\{X_{i,d}\}$.  Assume that the primary tree \(T\) for \(L\) is
normalized except at one or more vertices in the trunk of \(T\).
Let \(m\) be the number of carets in the trunk of \(T\).  Then
\(L\sim L'=C_{k_0,c_0}C_{k_1,c_1}\cdots C_{k_g,c_g}u'\) where
\(k_0<k_1<\cdots<k_g\), where \(u'\) is a word in the set $\{X_{i,d},\pi_s\}$, 
so that \(m\ge k_g+1\), and so that the complexity of the primary tree $T'$
of $L'$ is strictly less than the complexity of $T$.
%so that the smallest \(j\) so that \(i_j\ne k_j\) has \(i_j<k_j\).
}

\begin{proof}
Let \(\Lambda\) be the trunk of \(T\).  The interior
vertices of \(\Lambda\) are the interior, left edge vertices of
\(T\) and let these be \(a_0,a_1,\cdots, a_{m-1}\).  Let \(r\) be
the highest value with \(0\le r<m\) for which \(a_r\) is not
normalized.  Note that this is the lowest non-normalized interior
vertex of \(\Lambda\) and that, since $a_r$ is not normalized it is labelled $\ell \ne 1$
and must correspond to some $C_{i_j,\ell}$
and from Lemma \ref{thm:brin-lemma-4.15}, we have \(i_j=r\).

Moreover, since it is not normalized, $a_r$ must correspond to some hypercube $S_{i_j}$
which is fully divided across dimension $\ell$ and some other dimension $d$, with $1 \le d < \ell$.

By rewriting $L$ as $(t,s_{0,1}^k)$ (which we can do by Lemma \ref{thm:brin-lemma-4.15}) 
and applying Corollary \ref{thm:pull-up-fully-divided} 
to $t$, we can assume that the children of $a_r$, $v_1$ and $v_2$, are both labelled $d$.
We divide our work in two cases, $d=1$ and $d>1$. We observe that
the case $d=1$ is entirely analog to the proof of Theorem 4.22 in Brin \cite{B2}
while the case $d>1$ is slightly different.

In the case $d=1$, the left child $v_1$, which is in the trunk $\Lambda$, is labelled $1$. In the case that 
$j<n$ we observe that $i_{j+1} > r+1=i_j+1$, since the interior vertex of the trunk corresponding to $C_{i_{j+1},d_{j+1}}$
is not labelled $1$ 
%(A PICTURE HERE WOULD HELP SHOWING THAT THE LABEL OF $i_j$ is $\ell$, THE ONE 
% BELOW IT IS $1$, AND THE LABEL OF $C_{i_{j+1},d_{j+1}}$ IS SOME NUMBER 
% DIFFERENT FROM 1). 
Since the right child $v_2$ is an interior vertex not on the trunk, there must be a letter $X_{q,1}$
corresponding to it. By Lemma \ref{thm:brin-lemma-2.9} we can assume that $X_{q,1}$ occurs as the first letter of $u$,
that is $u=X_{q,1}u''$. Hence
\[
L=C_{i_0}\cdots C_{i_{j-1}}C_{i_j,\ell}C_{i_{j+1}} \cdots C_{i_g}X_{q,1}u''
\]
where we have omitted all the dimension subscripts of the baker's maps $C_{i,d}$ (except for one map) since they are not important
for the argument. The subword $C_{i_0}\cdots C_{i_j,\ell} \cdots C_{i_g}X_{q,1}$ is a trunk with a single caret labelled $1$ attached
at the caret $i_j$ of the trunk on its right child. By a careful observation of the right-left ordering it is evident that $q=i_j$.
By using relation (15) repeatedly on $L$ we can rewrite it as
\[
C_{i_0}\cdots C_{i_{j-1}}C_{i_j,\ell} X_{i_j,1}C_{i_{j+1}+1} \cdots C_{i_g+1}u'',
\]
since $i_0 < i_1 < \ldots < i_g$ and $i_{j+1} > i_j + 1$. 
Combining relations (15) and (16) on the product $C_{i_j,\ell} X_{i_j,1}$ we rewrite $L$ as
\[
C_{i_0}\cdots C_{i_{j-1}}(C_{i_j+1,\ell} X_{i_j,\ell}\pi_{i_j+1})C_{i_{j+1}+1} \cdots C_{i_g+1}u''.
\]
Now we apply (17) to commute $\pi_{i_j+1}$ back to the right without affecting the indices of the baker's maps.
This is possible since $i_{j+1}>i_j+1$ and therefore $i_{j+1}+1>i_j+2$. Now we apply (15) repeatedly
to the word
\[
C_{i_0}\cdots C_{i_{j-1}}C_{i_j+1,\ell} X_{i_j,\ell}C_{i_{j+1}+1} \cdots C_{i_g+1} \pi_{i_j+1} u''
\] 
to bring $X_{i_j,\ell}$ back to the right decreasing the indices of the the baker's maps by $1$ 
\[
C_{i_0}\cdots C_{i_{j-1}}C_{i_j+1,\ell} C_{i_{j+1}} \cdots C_{i_g} X_{i_j,\ell} \pi_{i_j+1} u''.
\]
By setting $u'=X_{i_j,\ell}\pi_{i_j+1} u''$ in the previous equation
and relabelling the indices with $k_i$'s, 
the word written in the previous equation has a primary tree $T'$ whose
complexity is strictly less than the complexity of $T$. The only thing
we still need to prove in this case is that $m \ge k_g + 1$. 
However, it has been observed above that $i_j = r < m - 1$ so $i_j + 2 \le m$. 
This gives the result in the case that $j = n$. If $j<n$, then $k_g =i_g$ and $m \ge i_g+1$ by Lemma
\ref{thm:brin-lemma-4.15}.

Now we study the case $1<d<\ell$. We observe that $a_r$ corresponds to $C_{i_j,\ell}$ and that
$v_1$ corresponds to $C_{i_k,d}$. By Lemma \ref{thm:brin-lemma-4.15}, we have $r+1=i_k$ which implies $i_k=i_j+1=i_{j+1}$.
In fact, if $i_j+1 < i_{j+1}$, there would be a vertex labelled $1$ on the trunk between the vertices $i_j$ and $i_{j+1}$
(and this is impossible since $d>1$). Let $X_{i_j,d}$ correspond to the right child $v_2$. Arguing as in the case $d=1$
we have 
\[
L=C_{i_0}\cdots C_{i_{j-1}}C_{i_j,\ell}C_{i_j+1,d}C_{i_{j+2}} \cdots C_{i_g}X_{q,d}u''
\]
and applying relation (15)
\[
L \sim C_{i_0}\cdots C_{i_{j-1}}(C_{i_j,\ell} C_{i_j+1,d} X_{i_j,d})C_{i_{j+2}} \cdots C_{i_g+1}u''.
\]
which can be rewritten as
\[
C_{i_0}\cdots C_{i_{j-1}} (C_{i_j,\ell} X_{i_j,d} C_{i_j+2,d})C_{i_{j+2}+1} \cdots C_{i_g+1}u''.
\]
By using the cross relation (18) on $C_{i_j,\ell} X_{i_j,d} C_{i_j+2,d}$ we read it as
\[
C_{i_0}\cdots C_{i_{j-1}} (C_{i_j,d} X_{i_j,\ell} C_{i_j+2,\ell} \pi_{i_j+1})C_{i_{j+2}+1} \cdots C_{i_g+1}u''
\]
Since $i_{j+2}>i_{j+1}$, then $i_{j+2}+1>i_{j+1}+1$, hence $\pi_{i_j+1}$ and the baker's maps to its right commute,
so the word becomes
\[
C_{i_0}\cdots C_{i_j,d} X_{i_j,\ell} C_{i_j+2,\ell} C_{i_{j+2}+1} \cdots C_{i_g+1}\pi_{i_j+1}u''.
\]
We apply (15) repeatedly and move $X_{i_j,\ell}$ back to the right to obtain
\[
L\sim C_{i_0}\cdots \underline{C_{i_j,d} C_{i_j+1,\ell}} C_{i_{j+2}} \cdots C_{i_g}X_{i_j,\ell}\pi_{i_j+1}u'',
\]
where the product $C_{i_j,d} C_{i_j+2,\ell}$ has been underlined to stress that the new trunk has the vertices
labelled $d$ and $\ell$ which are now switched. Thus the complexity of the tree has been lowered. In this second case, the new sequence $k_0 < \ldots < k_g$ 
is exactly equal to the initial one $i_0 < \ldots <i_g$. By the definition
of $m$ (given in Lemma \ref{thm:brin-lemma-4.15}) applied on the initial word $L$, 
we have that $m \ge i_g+1$ and so, since $k_g=i_g$, we are done.
\end{proof}

\remark{As observed in the proof, the case $d=1$ is equivalent to Theorem
4.22 in \cite{B2}, though the proof leads to a condition that is equivalent
to lowering the complexity.
When the index in some $C_{i_j,d}$ 
goes up by $1$, this corresponds to switching the vertices with labels
$d$ and $1$ in the primary tree and thus lowering the complexity 
by making more vertices normalized.}

\proposition{\label{thm:brin-lemma-4.19-and-4.20}Let 
$w$ be a word in the generating set 
\[
\{X_{i,d}, C_{i,d'}, \pi_i, \ov{\pi}_i,, X_{i,d}^{-1}, C_{i,d'}^{-1} 
\mid 1 \leq d \leq n, 2\le d' \le n, 
i \in \mathbb{N} \}. 
\]
Then $w \sim LMR$ as in Lemma \ref{thm:brin-lemma-4.8} and when expressed as elements of $\widehat{nV}$ we have $L = t s_{0,1}^{-p}$, $R^{-1} = y s_{0,1}^{-p}$, and $M = s_{0,1}^{p} u s_{0,1}^{-p}$ where $t$, $y$ are words in $\{s_{i,d} \mid 1 \leq d \leq n, i \in \mathbb{N}\}$, $u$ is a word in $\{ \sigma_j \mid 0 \leq j \leq p-1 \}$, and the lengths of $t$ and $y$ are both $p$. Further, we may assume the trees for $t$ and $y$ are normalized, and if $u$ can be reduced to the trivial word using relations (2) -- (4), then $M$ can be reduced to the trivial word using relations 
(13)--(17).}

\begin{proof} The proof of the first conclusion is exactly the same as the proof of lemma 4.19 of \cite{B3}. In order to assume the trees for $t$ and $y$ are normalized, we alternate applying Lemmas 
\ref{thm:brin-lemma-4.21} and \ref{thm:brin-lemma-4.22}. 
We have $L$ expressed as $(t, s_{0,1}^p)$, where $p$ is the length of $t$ and the number of carets in the trunk of the tree $T$ for $t$ is $m$. Setting $L = L'$ certainly gives that $L \sim L'$ and $m \geq k_g + 1$ by Lemma \ref{thm:brin-lemma-4.15}, 
so we have satisfied the hypotheses of Lemma 
\ref{thm:brin-lemma-4.21}. Therefore, $L \sim L_1z$ where $L_1$ expressed as $(t', s_{0,1}^p)$ where the trunk of the tree $T'$ for $t'$ has m carets. Since we set $L= L'$, we see that the trunks of $T$ and $T'$ are identical and the only way in which the two trees differ is that $T'$ is normalized off the trunk. Since $z$ is a word in $\{ \pi_i \}$, $z$ can be absorbed into $M$ without disrupting the assumptions on $M$, namely $M$ can still be written in the form $M = s_{0,1}^{p} u s_{0,1}^{-p}$ as above. We now replace $L$ with $L_1$ and proceed to use Lemma
\ref{thm:brin-lemma-4.22}.

Since the tree for $L$ is now normalized off the trunk, we satisfy the hypotheses of Lemma \ref{thm:brin-lemma-4.22} and write $L \sim L'$ where the tree for $L'$ has complexity lower than the tree for $L$ and $m \geq k_g +1$. Hence, we can now apply Lemma \ref{thm:brin-lemma-4.21} again and obtain $L\sim L_1z$ and let $z$ be absorbed into $M$. We apply this process over and over, decreasing the complexity of the tree associated to $L$ each time. Since there are only finitely many linearly ordered complexities, eventually this process will terminate, at which point the tree for $L$ will be normalized. We can apply the same procedure to the inverse of $LMR$ to normalize the tree for $R$. The last statement regarding $M$ follows immediately from Lemma 4.18 of \cite{B2}.

\end{proof}

\theorem{Let $w$ be a word in the generating set 
\[
\{X_{i,d}, C_{i,d'}, \pi_i, \ov{\pi}_i,, X_{i,d}^{-1}, C_{i,d'}^{-1} 
\mid 1 \leq d \leq n, 2\le d' \le n, 
i \in \mathbb{N} \}
\]
that represents the trivial element of $nV$. Then $w \sim 1$ using the relations in (1)--(18). Hence, we have a presentation for $nV$.}

\begin{proof}
Using the Proposition \ref{thm:brin-lemma-4.19-and-4.20}, we can assume 
\[
w \sim LMR = (t s_{0,1}^{-p}) (s_{0,1}^{p} u s_{0,1}^{-p}) (s_{0,1}^{p} y^{-1}) = tuy^{-1}
\]
where $t$, $y$ are words in $\{s_{i,d} \mid 1 \leq d \leq n, i \in \mathbb{N}\}$, $u$ is a word in $\{ \sigma_j \mid 0 \leq j \leq p-1 \}$, and the trees associated to $t$ and $y$ are normalized. By assumption, $tuy^{-1} = (tu, y)$ is the trivial element of $\widehat{nV}$ and so $tu$ and $y$ represent the same numbered patterns in $\Pi_n$. Furthermore, $t$ and $y$ must give the same unnumbered pattern, while $u$ enacts a permutation on the numbering. Since the forests for $t$ and $y$ are normalized and give the same pattern, the forests are identical with the same labeling by 
Lemma \ref{thm:brin-lemma-2.10}. The numbering on the leaves for both forests follows the left-right ordering, hence $t$ and $y$ give the same numbered patterns, which implies that $u$ enacts the trivial permutation and $M \sim 1$ by Proposition \ref{thm:brin-lemma-4.19-and-4.20}.

We now wish to show that $L \sim R^{-1}$. By Lemma \ref{thm:brin-lemma-4.8}, we have

\begin{itemize}
\medskip
\item $L=C_{i_0, d_0}C_{i_1, d_1}\dots C_{i_g, d_g}q$ 
\medskip
\item $R^{-1}=C_{j_0, d'_0}C_{j_1, d'_1}\dots C_{j_m, d'_m}q'$ 
\medskip
\end{itemize}

Since we know that the trunks of the trees corresponding to $L$ and $R^{-1}$ are identical with the same labeling, the sequences $(i_0, i_1, \dots, i_g)$ and $(j_0, j_1, \dots, j_m)$ are identical and $d_k = d'_k$ for each $k \in \{ 0, 1, \dots, n=m \}$. Hence, the subwords $C_{i_0, d_0}C_{i_1, d_1}\dots C_{i_g, d_g}$ and $C_{j_0, d'_0}C_{j_1, d'_1}\dots C_{j_m, d'_m}$ are the same and it remains to show that $q \sim q'$. This follows from Lemma \ref{thm:brin-lemma-2.8} 
and the homomorphism from $\widehat{nV}$ to $nV$ as before.

 \end{proof}

\section{\label{sec:finite-presentations}Finite Presentations}

\subsection{\label{sec:finite-presentation}Finite Presentation for $\widehat{nV}$}

We now give a finite presentation for $\widehat{nV}$, using analogous arguments found in \cite{B2} to show that the full set of relations is the result of only finitely many of them. First, recall our generating set is $\{s_{i,d}, \sigma_i \mid i \in \mathbb{N}, 1 \leq d \leq n \}$. When $i < j$, relations $(\M1)$ and $(\M5a)$ give $s_{i,1}^{-1}x_j s_{i,1} = x_{j+1}$ where $x_j = s_{j,d}$ 
(for some $d$) or $\sigma_j$. Hence, we can use $s_{i,d} = s_{0,1}^{1-i}s_{1,d}s_{0,1}^{i-1}$ and $\sigma_i = s_{0,1}^{1-i}\sigma_1s_{0,1}^{i-1}$ as definitions for $i \geq 2$. Therefore, $\widehat{nV}$ is generated by $\{s_{i,d}, \sigma_i \mid i \in \{0,1\} , 1 \leq d \leq n \}$, which gives a generating set of size $2n + 2$ for each $n$.

We treat relations $(\M1)$ through $(\M6)$ in the same way as they are treated in \cite{B2}. Relations involving only one parameter, such as $(\M2)$, $(\M4)$, and $(\M6)$, 
are obtained for $i \geq 2$ by setting $i = 1$ and conjugating by powers of $s_{0,1}$, therefore the only necessary relations to include are when $i = 0$ and $i = 1$. As before, $(\M2)$ and $(\M4)$ follow from: $\sigma_0^2 = 1$, $\sigma_1^2 = 1$, $\sigma_0 \sigma_1 \sigma_0 = \sigma_1 \sigma_0 \sigma_1$, and $\sigma_1 \sigma_2 \sigma_1 = \sigma_2 \sigma_1 \sigma_2$, or 4 relations for each $n$. Relation (7) follows from 2 relations for each pair of distinct dimensions, giving $2 {n \choose 2} = n(n-1)$ relations for each $n$.

Relation $(\M3)$ is treated the same way as in \cite{B2} for each $n$. Hence, for all $i, j$, $(\M3)$ follows from the 4 relations: 
$\sigma_0 \sigma_2 = \sigma_2 \sigma_0$, $\sigma_0 \sigma_3 = \sigma_3 \sigma_0$, $\sigma_1 \sigma_3 = \sigma_3 \sigma_1$, $\sigma_1 \sigma_4 = \sigma_4 \sigma_1$.

For relation $(\M1)$, which can be rewritten as $s_{i,d}^{-1}s_{i + k,d'}s_{i,d}=s_{i+k+1,d'}$ for $k>0$, we have two cases: the case where $d = 1$ and the case where $d \neq 1$. If $d=1$, then the case $i = 0$ follows by definition, and by the same induction argument used in \cite{B2} implies that the relation for all $i,k$ follows from the cases where $i = 1$ and $k = 1, 2$, hence we need only 2 relations per dimension. If $d \neq 1$, we do not get the case $i = 0$ by definition and we must include $i =0,1$ and $k =1, 2$, i.e. 4 relations per each pair of dimensions. There are $n-1$ choices for $d$, as $d \neq 1$, and $n$ choices for $d'$, so this case yields $4n(n-1)$ relations. Hence, in total $(\M1)$ 
can be obtained for all $i, k$ by $2n + 4n(n-1) = 4n^2 -2n$ relations.

For relation $(\M5b)$, $\sigma_{i} s_{i,d}=s_{i+1,d}\sigma_i\sigma_{i+1}$, there is only a single parameter to deal with, hence the relation for $i \geq 2$ can be obtained from the cases where $i=0,1$ by conjugating by $s_{0,1}$ as before. Relation $(\M5c)$
is actually equivalent to $(\M5b)$, hence for each $n$ we only need $2n$ relations for $(\M5b)$, $(\M5c)$. We treat $(\M5a)$ $\sigma_{i+k}s_{i,d}=s_{i,d}\sigma_{i+k+1}$ for $k>0$ the same way as for $(\M1)$, hence 2 relations are required for $d = 1$ and 4 for $d \neq 1$ for a total of $4n - 2$ relations. And lastly, $(\M5d)$ $\sigma_{i}s_{i+k,d}=s_{i+k,d}\sigma_i$ can be obtained in the same way as the second case of $(\M1)$ where the relation for all $i, k$ is obtained by $i = 0,1, k = 2,3$, i.e. $4n$ relations.

Thus, we have proven the following:

\theorem{The group $\widehat{nV}$ is presented by the $2n + 2$ generators $\{s_{i,d}, \sigma_i \mid i \in \{0,1\} , 1 \leq d \leq n \}$ and the $5n^2 + 7n +6$ relations given below:}

\[
\begin{array}{cccc}
(\M1) & & s_{1,1}^{-1}s_{1 + k,d'}s_{1,1}=s_{2 + k,d'} & k= 1,2 \medskip \\

    & & s_{i,d}^{-1}s_{i + k,d'}s_{i,d}=s_{i+k+1,d'} &  i = 0,1, k=1,2; 2\leq d \leq n \medskip \\

(\M2) & & {\sigma_i}^2 = 1 & i = 0,1 \medskip \\

(\M3) & & \sigma_i \sigma_{i+k} = \sigma_{i+k} \sigma_i & i = 0,1, k = 2, 3 \medskip \\

(\M4) & & \sigma_i\sigma_{i+1}\sigma_i=\sigma_{i+1}\sigma_i\sigma_{i+1}& i = 0,1 \medskip \\

(\M5a) & & \sigma_{k+1}s_{1,1}=s_{1,1}\sigma_{k+2} & k = 1,2 \medskip \\

        & & \sigma_{i+k}s_{i,d}=s_{i,d}\sigma_{i+k+1} & i = 0,1, k = 1, 2;  2\leq d \leq n \medskip \\

(\M5b / \M5c) & & \sigma_{i} s_{i,d}=s_{i+1,d}\sigma_i\sigma_{i+1} & i= 0,1 \medskip \\

(\M5d) & & \sigma_{i}s_{i+k,d}=s_{i+k,d}\sigma_i & i = 0,1, k =2,3 \medskip \\

(\M6) & & s_{i,d}s_{i+1,d'}s_{i,d'}=s_{i,d'}s_{i+1,d}s_{i,d}\sigma_{i+1} & i = 0,1, d \neq d' \\

\end{array}
\]

\subsection{Finite Presentation for $nV$} We can use the relations in $nV$ to write
\[\begin{split}
X_{i,d}&=X_{0,1}^{1-i}X_{1,d}X_{0,1}^{i-1}, \\
\pi_i&=X_{0,1}^{1-i} \pi_1X_{0,1}^{i-1}, \\
\overline{\pi}_i&=X_{0,1}^{1-i}\overline{\pi}_1X_{0,1}^{i-1}
\end{split}
\] 
for \(i\ge2\) and $1 \leq d \leq n$. We can also use the relations for $nV$ as in Proposition 6.2 of \cite{B1} to write
\[
\begin{split}
C_{m,d}&=(\overline{\pi}_m X_{m,d} \overline{\pi}_{m+1} \pi_m)(X_{m,d} \pi_{m+1} X_{m,1}^{-1}) 
\end{split}
\] 
for $m \geq 0$ and $2 \le d \le n$, which we use as a definition. Hence, the $C_{m,d}$ are not needed to generate $nV$. 

The homomorphism $\widehat{nV} \rightarrow nV$ given by $s_{i,d} \mapsto X_{i,d}$ and $\sigma_i \mapsto \pi_i$ implies that the work done for the relations for $\widehat{nV}$ carries over to relations (1)--(4), (7)--(9), and (12)
(see Lemma \ref{thm:brin-lemma-4.1}). 
Relations (10)-(11) and (13)-(6) are exactly the same as those from $2V$ and can be treated as in \cite{B2}, contributing a total of 10 relations to our finite set. 

Relation (5) can be treated in a manner similar to $(\M1)$ from $\widehat{nV}$, where 2 relations are needed for dimension 1 and 4 for all others, contributing a total of $4(n-1) + 2$ relations. Relations (14) and (16) include only one parameter and hence can be obtained from the cases where $i = 0,1$ as before, contributing $2(n-1)$ relations apiece. And (17) requires 4 relations for each $d \neq 1$, hence adding an additional $4(n-1)$ relations.

For relation (15), we have two cases: for $d'=1$, all cases follow from when $i =0,1$, giving us $2(n-1)$ relations since $2 \leq d \leq n$. For $d' \neq 1$, 4 relations are required for each pair $d, d' \in \{2, . . . , n\}$, contributing $4(n-1)(n-1)$ relations. And lastly, since (18) involves only one parameter in the first component, we only need 2 relations for each $1 < d' < d \le n$, the number of such pairs being 
$\frac{(n-1)(n-2)}{2}$. 

We now have the following:

\theorem{\label{thm:finite-presentation}
The group $nV$ is presented by the $2n + 4$ generators $\{X_{i,d}, \pi_i, \overline{\pi}_i \mid i \in \{0,1\}, 1 \leq d \leq n\}$, the $5n^2 + 7n +6$ relations obtained from the homomorphism $\widehat{nV} \rightarrow nV$, and the additional $5n^2 + 3n + 4$ relations given below for a total of $10n^2 + 10n + 10$ relations.}

\[
\begin{array}{cccc}

(5) & & \overline{\pi}_{k+1} X_{1,1} = X_{1,1} \overline{\pi}_{k+2} & k = 1,2 \medskip \\

      & & \overline{\pi}_{m+k} X_{m,d} = X_{m,d} \overline{\pi}_{m+k+1} & m=0,1, k =1,2, 2 \leq d \leq n \medskip \\

(10) & & \overline{\pi}_{m+k} \pi_m = \pi_m \overline{\pi}_{m+k} & m =0,1, k =2,3 \medskip \\

(11) & & \pi_m \overline{\pi}_{m+1} \pi_m = \overline{\pi}_{m+1} \pi_m \overline{\pi}_{m+1}  & m =0,1 \medskip \\

(13) & & \overline{\pi}_m^2 =1 & m =0,1 \medskip \\

(6) & &  \overline{\pi}_m X_{m,1} = \pi_m \overline{\pi}_{m+1} & m =0,1 \medskip \\

(14) & &  \overline{\pi}_m X_{m,d} = C_{m+1,d} \pi_m \overline{\pi}_{m+1}  & m= 0,1, d \ne 1 \medskip \\

(15) & & C_{k+1,d} X_{1,1} = X_{1,1} C_{k+2,d} & k=1,2  \medskip \\

        & & C_{m+k,d} X_{m,d'} = X_{m,d'} C_{m+k+1,d} & m=0,1, k=1,2; 2 \le d, d' \le n \medskip \\

(16) & & C_{m,d} X_{m,1} = X_{m,d} C_{m+2,d} \pi_{m+1} & m= 0,1; 2 \le d \le n \medskip \\

(17) & & \pi_{m} C_{m+k,d} = C_{m+k,d} \pi_{m} & m=0,1, k =2,3 \medskip \\

(18) & & C_{m,d} X_{m,d'} C_{m+2,d'}= C_{m,d'} X_{m,d} C_{m+2,d} \pi_{m+1} & m= 0,1; 1 < d' < d \le n \medskip \\
\end{array}
\]

\remark{Since $\omega V$ is an ascending union of the $nV$'s, a word 
$w \in \{X_{i,d}, \pi_i, \overline{\pi}_i \mid i \in \{0,1\}, d \in \mathbb{N} \}$ such that $w =_{\omega V} 1$
must be contained in some $nV$ (for some $n \in \mathbb{N}$) and so we can use the same ideas
and the relations inside $nV$ to transform $w$ into the empty word. Therefore, the following result is
an immediate consequence of Theorem \ref{thm:finite-presentation}.}

\corollary{The group $\omega V$ is generated by the set 
$\{X_{i,d}, \pi_i, \overline{\pi}_i \mid i \in \{0,1\}, d \in \mathbb{N} \}$ and satisfies the family of relations in Theorem
\ref{thm:finite-presentation} with the only exception that the parameters $d,d' \in \mathbb{N}$.}

\section{\label{sec:simplicity}Simplicity of $nV$ and $\omega V$}

Brin proved in \cite{B3} that the groups $nV$ and $\omega V$ are simple by
showing that the baker's map is a product of transpositions and following the outline of
an existing proof that $V$ is simple. 

We reprove Brin's simplicity result verify that Brin's original 
proof that $2V$ is simple (Theorem 7.2 in \cite{B1})
generalizes using the generators and the relations that have been found.

\theorem{\label{thm:group-equal-commutators}The 
groups $nV$, $n \le \omega$, equal their commutator subgroups.}

\begin{proof}
The goal is to show that the generators $X_{m,i}, \pi_m, \ov{\pi}_m$ are products
of commutators. We write $f \simeq g$ to mean that $f=g$ modulo the commutator subgroup. We also observe that the arguments below are independent of the dimension
$i$.

From relation (1) we see that 
$X_{q,i}^{-1}X_{0,1}^{-1}X_{q,i}X_{0,1}=X_{q,i}^{-1} X_{q+1,i}$ for $q \ge 1$
and so $X_{q+1,i} \simeq X_{q,i}$. Therefore $X_{q,i} \simeq X_{1,i}$, for $q \ge 1$.
Using relation (2) and arguing similarly, we see that $\pi_q \simeq \pi_1$, for
$q \ge 1$.

From relation (3) we see that 
$\pi_0 X_{0,i} \pi_0^{-1} X_{0,i}^{-1}= X_{1,i}\pi_1 X_{0,i}^{-1}$
so that $X_{0,i} \simeq X_{1,i}\pi_1$. 
Also, by relation (3), $X_{2,i} \simeq X_{1,i}$ and 
the fact that
$\pi_2 \simeq \pi_1$, we see $\pi_1 X_{1,i} = X_{2,i}\pi_1 \pi_2 \simeq X_{1,i}\pi_1 
\pi_1= X_{1,i}$. Therefore $\pi_1 \simeq 1$ and so $X_{0,i} \simeq X_{1,i}$.

Relation (9) and $\pi_1 \simeq 1$ give that  $\pi_0^2 \simeq \pi_0\pi_1\pi_0=\pi_1\pi_0\pi_1 \simeq \pi_0$
which implies $\pi_0\simeq1$.

By relation (6) and the fact that 
$\pi_1 \simeq 1$ and $\ov{\pi}_1 \simeq \ov{\pi}_0$ we get
$\ov{\pi}_1 X_{1,1} = \pi_1 \ov{\pi}_2 \simeq \ov{\pi}_1$. Hence
$X_{0,1}\simeq X_{1,1} \simeq 1$. 

Now, relation (6) and $X_{0,1} \simeq 1$ give that
$\ov{\pi}_0 \simeq \ov{\pi}_0 X_{0,1} = \ov{\pi}_1$.
Relation (11) and $\pi_0 \simeq 1$ lead to
$\ov{\pi}_1 \simeq \pi_0 \ov{\pi}_1 \pi_0 = \ov{\pi}_1 \pi_0 \ov{\pi}_1 \simeq \ov{\pi}_1^2$.
Therefore $\ov{\pi}_0 \simeq \ov{\pi}_1 \simeq 1$.

Finally, by relation (7) and $X_{0,1} \simeq X_{1,1} \simeq 1 \simeq \pi_1$
we get 
$X_{1,i} X_{0,i} \simeq X_{0,1} X_{1,i} X_{0,i} = X_{0,i}  X_{1,1} X_{0,1} \pi_1
\simeq X_{0,i}$ which implies $X_{0,i} \simeq X_{1,i} \simeq 1$.
We have thus proved that all the generators of $nV$
are in the commutator subgroup. The case of $\omega V$ is identical: each generator
lies in some $nV$ and can be written as a product of commutators within that
subgroup.
\end{proof}

From Section 3.1 in \cite{B1} (which generalizes to $nV$ and $\omega V$
as observed by Brin in \cite{B2} and \cite{B3}) the commutator subgroup of $nV$ and $\omega V$ are simple, therefore Theorem \ref{thm:group-equal-commutators}
implies the following result.

\theorem{The groups $nV$, $n \le \omega$, are simple.}

\section{An alternative 
generating set \label{sec:alternative-generating-set}}

Observe that, 
for any $n \in \mathbb{N}$, we have $(n-1)V \times V \le nV$.
It can be shown that another generating set for $nV$ is given by taking
a generating set for $(n-1)V \times V$ and adding an involution which
swaps two disjoint subcubes of $[0,1]^n$, one of which has the origin as one of its
vertices and the other one which contains the vertex $(1, \dots, 1)$.
This second generating set has the advantage of taking the generators of
$(n-1)V$ and adding only the generators of $V$ plus another one. This leads
to a smaller generating set which was suggested to us by Collin Bleak. 
It seems feasible that a good set of relations exist for this alternative generating set.

%\section{An alternative finite generating set}
%
%Brief discussion of Collin's finite generating set, Ruskuc result to keep the number of %generators bounded independently
%of the dimension of $nV$.


\begin{thebibliography}{1}

\bibitem{BL} C. Bleak and D. Lanoue, {\it A family of non-isomorphism results}, Geom. Dedicata $\mathbf{146(1)}$ (2010), 21--26.
Arxiv preprint: \url{http://arxiv.org/abs/0807.4955}

\bibitem{B1} M. G. Brin, {\it Higher Dimensional Thompson Groups}, Geom. Dedicata $\mathbf{108}$ (2004) 163-192.
Arxiv preprint: \url{http://arxiv.org/abs/math/0406046}

\bibitem{B2} M. G. Brin, {\it Presentations of higher dimensional Thompson Groups}, J. Algebra $\mathbf{284}$ (2005) no. 2, 520--558. 
Arxiv preprint: \url{http://arxiv.org/abs/math/0501082}

\bibitem{B3} M. G. Brin, {\it On the baker's map and the simplicity of the higher dimensional Thompson groups $nV$}, Publ. Mat. $\mathbf{54}$ (2010) 433-439.
Arxiv preprint: \url{http://arxiv.org/abs/0904.2624}

\bibitem{CFP} J. W. Cannon, W. J. Floyd, and W. R. Parry, {\it Introductory Notes on Richard Thompson's groups}, Enseign. Math. (2) $\mathbf{42}$ (1996) 215-256 \

\bibitem{GKKL} R. M. Guralnick, W. M. Kantor, M. Kassabov, A. Lubotzky,
Presentations of finite simple groups: a computational approach, \emph{to appear},
\url{http://arxiv.org/abs/0804.1396}

\bibitem{KMN} D. H Kochloukova, C. Martinez-Perez, and B. E. A. Nucinkis, {\it Cohomological finiteness properties of the Brin-Thompson-Higman groups 2V and 3V}, ArXiv preprint: \url{http://arxiv.org/abs/1009.4600}.



\end{thebibliography}
\end{document}